\documentclass{conm-p-l}
\usepackage{amssymb, hyperref}
\usepackage{xcolor}

\newtheorem{theorem}{Theorem}[section]

\newtheorem{problem}{Problem}[section]
\newtheorem{proposition}{Proposition}[section]
\newtheorem{remark}{Remark}[section]
\usepackage{subfig, graphicx}
\usepackage{algorithm}
\usepackage{algpseudocode}

\theoremstyle{definition}

\theoremstyle{remark}

\numberwithin{equation}{section}

\newcommand{\x}{{\bf x}}
\newcommand{\R}{\mathbb{R}}
\newcommand{\C}{\mathbb{C}}
\newcommand{\ii}{{\rm i}}
\newcommand{\ik}{\ii k}

\begin{document}

\title[The QRM for an ISP]{ Reconstructing a space-dependent source term via the quasi-reversibility method}

\author[L. Nguyen]{Loc H. Nguyen}
\address{Department of Mathematics and Statistics, University of North Carolina at Charlotte, Charlotte
, North Carolina, USA, 28223}
\email{loc.nguyen@uncc.edu}
\thanks{The work of LHN was supported by US Army Research Laboratory and US Army Research Office grant W911NF-19-1-0044 and was supported, in part, by funds provided by the Faculty Research Grant program at UNC Charlotte.}

\author[H. Vu]{Huong T.T. Vu}
\address{Department of Information Technology, University of Finance-Marketing, Ho Chi Minh City, Vietnam}
\email{vtthuong@ufm.edu.vn}

\subjclass{Primary 35R30; Secondary: 78A46}
\date{}


\keywords{Inverse source problem, truncated Fourier series, approximation, Carleman estimate, convergence}

\begin{abstract}
The aim of this paper is to solve an important inverse source problem which arises from the well-known inverse scattering problem.
We propose to truncate the Fourier series of the solution to the governing equation with respect to a special basis of $L^2$. By this, we obtain a system of linear elliptic equations.
Solutions to this system are the Fourier coefficients of the solution to the governing equation.
After computing these Fourier coefficients, we can directly find the desired source function.
Numerical examples are presented.
\end{abstract}

\maketitle

%
.

\section{Introduction}
Let $d \geq 1$ be the spatial dimension.
Fix a wavenumber $k > 0$.
Let $[\underline \theta, \overline \theta]$ be an interval of angles. 
Let $u_0 = u_0(\x, \theta)$, $(\x, \theta) \in \R^d \times [\underline \theta, \overline \theta]$, be the incident  plane wave that illuminates the medium at an angle $\theta$.
This choice of $u_0$ as the incident wave arises from the well-known  inverse scattering problem, stated later.
Let $u(\x, \theta)$, $(\x, \theta) \in \R^d \times [\underline \theta, \overline \theta]$, be the wave function that  is governed by the Helmholtz equation and the Sommerfeld radiation condition
\begin{equation}
	\left\{
		\begin{array}{ll}
			\Delta u(\x, \theta) + k^2 c_0(\x)u(\x, \theta) = -k^2 p(\x) u_0(\x, \theta) &\x \in \R^d,\\
			\partial_{|\x|} u(\x, \theta) - \ik u(\x, \theta) = o(|\x|^{\frac{1-d}{2}})
			&|\x| \to \infty.
		\end{array}
	\right.
	\label{Heqn}
\end{equation}
Here,  $c_0: \R^d \to [1, \infty)$ and $u_0: \R^d \times [\underline \theta, \overline \theta] \to \C$ are given functions.
The function $c_0$ represents the background dielectric constant of the medium.

Let $\Omega$ be an open and bounded domain of $\R^d$ with a smooth boundary.	Assume that the source function $p$ is compactly contained in $\Omega$.
We are interested in the following problem.
\begin{problem}[Inverse source problem]
	Determine the source function $p(\x)$ for all $\x \in \Omega$ from the following boundary measurements
	\begin{equation}
		f(\x, \theta) = u(\x, \theta)
		\quad 
		\mbox{and}
		\quad
		g(\x, \theta) = \partial_{\nu} u(\x, \theta)
		\label{1.2}
	\end{equation}
	for all $(\x, \theta) \in \partial \Omega \times [\underline \theta, \overline \theta].$
	Here, $\nu(\x)$ is the unit normal vector to $\partial \Omega$ at $\x$.
\label{isp}
\end{problem}

This inverse source problem is the linearization of the nonlinear inverse scattering problem. Since the inverse scattering problem has many real-world applications; mostly in exploring some inaccessible regions from external measurement;
for e.g, bio-medical imaging, nano-sciences, security, seismic exploration.
Therefore, Problem \ref{isp} is significant in those fields. 
The solver for Problem \ref{isp} proposed in this paper consists of three steps.
\begin{enumerate}
    \item \label{step1} We derive a system of linear elliptic equations. Solution of this system is a vector involving the first $N$ Fourier coefficients of the function $w = -u/(k^2 u_0)$, $ N \in \mathbb{N}$, with respect to a special basis of $L^2$ introduced in \cite{Klibanov:jiip2017}.
    \item \label{step2}  We apply the quasi-reversibility method developed in \cite{NguyenLiKlibanov:2019} to solve the system obtained in Step \ref{step1}. 
    \item Directly compute the desired source from the solution obtained in Step \ref{step2}.
\end{enumerate}

The inverse source problem and some of its versions were studied intensively. We cite to \cite{BaoLinTriki:jde2010, BaoLinTriki:CRM2011, BaoLinTriki:cm2011, EntekhabiIsakov:ip2018, IsakovLu:SIAM2018, IsakovLu:ipi2018} for the uniqueness, stability and numerical methods to solve inverse source problems for the case when the medium is homogeneous; i.e. $c_0 = \mbox{constant}$. 
The numerical reconstruction methods in those publications are based on the least squares optimization method.
Good quality  reconstructions due to those approaches are achieved only when the wave number $k$ is large. 
We observe that when $k$ is large, the data is  very sensitive with noise. This is because of the high oscillation of the data. 
Unlike this, the reconstructive method  in this paper and the method in \cite{NguyenLiKlibanov:2019} do not require the data at high frequency. Therefore, we can reconstruct the source with reasonable value of $k$. The difficulty about the noise is overcome.

As mentioned in step \ref{step2} above, we only compute the first $N$ Fourier coefficients of the function $w = -u/(k^2 u_0)$. 
That means we only solve Problem \ref{isp} in an ``Galerkin" approximation context. 
Rigorously verifying the convergence of this approximation as $N \to \infty$ is extremely challenging. We assume that this approximation is valid.
In contrast, this numerical approach is very effective  for many kinds of inverse problems; see e.g. \cite{Nguyen:CAMWA2020, KlibanovNguyen:ip2019, NguyenKlibanov:ip2022, Nguyens:jiip2020, LeNguyenNguyenPowell:JOSC2021}. This is the reason we employ this technique again in this paper.

The paper is organized as follows. In Section \ref{motivation}, we provide in details the formulation of Problem \ref{isp}.
In Section \ref{sec_Ana}, we establish the approximation context and derive a system of PDEs that plays an important role in our algorithm.
In Section \ref{s5}, we recall the quasi-reversibility method to solve the system obtained in Section \ref{sec_Ana}. We present numerical study in Section \ref{sec_num}. Section \ref{sec_concl} is for concluding remarks.

\section{The significance of Problem \ref{isp}} \label{motivation}

Problem \ref{isp} is the linearization  of the nonlinear inverse scattering problem, which has many real-world applications. We will list some of the important applications later.
Let $c: \R^d \to [1, \infty)$ be the spatially distributed dielectric constant of the medium.
For each angle $\theta \in [\underline \theta, \overline \theta],$   define 
\[
	\xi(\theta) = (\cos \theta, \sin \theta, 0, \dots, 0) \in \mathbb{S}^{d-1} 
\] that represents the direction of the angle $\theta$ in the $x_1 x_2$-plane.
Here, $\mathbb{S}^{d-1} = \{\x \in \R^d, |\x| = 1\}$ is the unit sphere in $\R^d.$
 We use the plane wave, so-called the incident plane wave, of the form
\begin{equation}
	v_0(\x, \theta) = e^{\ik \x \cdot \xi(\theta)} 
	\quad
	\mbox{for all } (\x, \theta) \in \R^d \times [\underline \theta, \overline \theta], k > 0,
	\label{incwave}
\end{equation}
to illuminate the medium.
The incident wave propagates in space and scatters.
 The resulting total wave $v(\x, \theta)$ is governed by the Helmholtz equation 
\begin{equation}
	\Delta v(\x, \theta) + k^2 c(\x) v(\x, \theta) = 0 
	\quad \mbox{for all } \x \in \R^d, \theta \in [\underline \theta, \overline \theta]
	\label{2.1}
\end{equation}
and the Sommerfeld radiation condition
\[
	\partial_{|\x|} v_{\rm sc}(\x, \theta) - \ik v_{\rm sc}(\x, \theta) =  o(|\x|^{\frac{1 - d}{2}})
\]
where 
\begin{equation}
	v_{\rm sc}(\x, \theta) = v(\x, \theta) - v_0(\x, \theta)
	\quad
	\mbox{for all } (\x, \theta) \in \R^d \times [\underline \theta, \overline \theta]
\end{equation}
is the scattering wave.

	By \eqref{incwave}, the angle between the $x_1$ axis and the propagation of incident wave in the $x_1 x_2$ plane is $\theta$. Roughly speaking, $\theta$ can be considered as the angle of view. Since we consider the case of complete angles of view, we  choose $[\underline \theta, \overline \theta] = [0, 2\pi].$ The study of Problem \ref{isp} in the case of partial angle of views will be studied  later.

The inverse scattering problem is formulated as follows.
\begin{problem}[The inverse scattering problem]
Let $\Omega$ be a bounded and open domain of $\R^d$ with smooth boundary.
Given the boundary measurement $v_{\rm sc}(\x, \theta)$ and $\partial_{\nu} v_{\rm sc}(\x, \theta)$ for all $(\x, \theta) \in \partial \Omega \times [\underline \theta, \overline \theta],$
determine the spatially distributed dielectric constant $c(\x)$ for all $\x \in \Omega$, 
\label{cip}
\end{problem}

{\it Why is the inverse scattering problem interesting and significant?}
According to  the formulation of Problem \ref{cip},  we want to compute the information of the spatially distributed dielectric constant of the medium from the external measurement of the wave field.
 The knowledge of  the reconstructed spatially distributed dielectric constant of a medium directly provides significant information about unknown objects inside that medium; for e.g, position, shape, size and physical properties. Typical examples of those objects are anti-personnel explosive devices buried under the ground, cancerous tumors inside living tissues, and nano structures. 
Therefore, solving the inverse scattering problem has important applications in  bio-medical imaging, nondestructive testing, radar, security, optical physics, seismic exploration, nano science.
Hence, the inverse scattering problem and the related ones have been intensively studied.
We list here several approaches: the imaging techniques based on sampling and the factorization methods \cite{ AmmariChowZou:sjap2016, AmmariKang:lnim2004,  Tan1:ip2012, Burge2005,   Colto1996, HarrisLiem:SIAM2020, Kirsc1998, LiLiuZou:smms2014, LiLiuWang:jcp2014, Liem:jiiptoappear, Soumekh:SAR}, the methods based on optimization \cite{Bakushinsii:springer2004, Tan2:ip2012, Chavent:springer2009, Engl:Kluwer1996, Gonch2013, Tihkonov:kapg1995}, the methods based on Born series \cite{Bleistein:ap1984, Chew:vnr1990, Devaney:cup2012, Kirsch:aa2017, Langenberg:1987, Moskow:Ip2008}, the method based on linearization \cite{BaoLi:ip2005, BaoLi:SIAM2005, Bao:ip2015, Chen:ip1997}, and the convexification method \cite{VoKlibanovNguyen:IP2020, Khoaelal:IPSE2021, KhoaKlibanovLoc:SIAMImaging2020}.
See \cite{ColtonKress:2013} for a more complete list.

The current paper contributes to the field by solving the linearization of the inverse scattering problem.
Assume that the dielectric constant $c$ is a perturbation of a background function $c_0$.
For simplicity, in this section, we choose $c_0$ to be identically 1. 
That means, the function $c$ is of the form
\begin{equation}
	c(\x) = 1 + \eta p(\x)
	\quad
	\mbox{for all } \x \in \R^d
	\label{2.3}
\end{equation}
where $0 < \eta \ll 1$  and $p$ is a function that indicate the unknown inclusion.
Therefore, the wave function $v$ satisfies
\begin{equation}
	\Delta v(\x, \theta) + k^2 (1 + \eta p(\x))v(\x, \theta) = 0
	\quad 
	\mbox{for all }
	(\x, \theta) \in \R^d \times [\underline \theta, \overline \theta].
	\label{2.4}
\end{equation}
It is not hard to verify that
\begin{equation}
	\Delta v_0(\x, \theta) + k^2 v_0(\x, \theta) = 0
	\quad 
	\mbox{for all }
	(\x, \theta) \in \R^d \times [\underline \theta, \overline \theta].
	\label{2.5}
\end{equation}
Let $u = v_{\rm sc}/\eta$. It follows from \eqref{2.4} and \eqref{2.5} that
\begin{equation}
	\Delta u(\x, \theta) + k^2 u(\x, \theta) = -k^2  p(\x) v(\x, \theta) 
	\quad
	\mbox{for all }
	(\x, \theta) \in \R^d \times [\underline \theta, \overline \theta].
	\label{2.6}
\end{equation}
It is obvious that $\lim_{\eta \to 0} v(\x, \theta) = v_0(\x, \theta)$. Thus,
since $\eta$ is a small number, equation \eqref{2.6} can be approximated as
\begin{equation}
	\Delta u(\x, \theta) + k^2 u(\x, \theta) = -k^2 p(\x) v_0(\x, \theta) 
	\quad
	\mbox{for all }
	(\x, \theta) \in \R^d \times [\underline \theta, \overline \theta].
\end{equation}
By choosing $u_0$ as the incident wave $v_0$, we derive Problem \ref{isp}.
While the argument above only give an example for the motivation to solve Problem \ref{isp}, our method to solve inverse problem is not limited in the context of the inverse scattering problem. That means, the proposed method in this paper can be applied in the general case when $u_0$ is not necessary to be the incident wave $v_0$.


\section{An approximation context} \label{sec_Ana}

Define 
\begin{equation}
	w(\x, \theta) = -\frac{u(\x, \theta)}{k^2 u_0(\x, \theta)}
	\quad 
	\mbox{or }
	\quad
	u(\x, \theta) = -k^2w(\x, \theta) u_0(\x, \theta)
	\label{w}
\end{equation}
for all $ (\x, k) \in \Omega \times [\underline \theta, \overline \theta].$
Since $\Delta u(\x, \theta) = -k^2 \Delta [w(\x, \theta) u_0(\x, \theta)]$, we have for all $(\x, \theta) \in \Omega \times [\underline \theta, \overline \theta]$
\begin{equation}
	\Delta u(\x, \theta) 
	= -k^2\big[w(\x, \theta) \Delta u_0(\x, \theta) + u_0(\x, \theta) \Delta w(\x, \theta) + 2\nabla w(\x, \theta) \cdot \nabla u_0(\x, \theta)\big].
	\label{3.1}
\end{equation}
On the other hand,
it follows from the Helmholtz equation in \eqref{Heqn} and the second identity in \eqref{w} that 
\begin{align*}
	\Delta u(\x, \theta) &= - k^2\big[ c_0(\x) u(\x, \theta) +  p(\x)u_0(\x, \theta)\big]\\
	&= - k^2 \big[-k^2 c_0(\x) w(\x, \theta) u_0(\x, \theta) + p(\x)u_0(\x, \theta)\big]
\end{align*}
for all $(\x, \theta) \in \Omega \times [\underline \theta, \overline \theta].$
This and \eqref{3.1} imply
\begin{multline*}
w(\x, \theta) \Delta u_0(\x, \theta) + u_0(\x, \theta) \Delta w(\x, \theta) + 2\nabla w(\x, \theta) \cdot \nabla u_0(\x, \theta)
\\
=
	- k^2 c_0(\x) w(\x, \theta) u_0(\x, \theta) + p(\x)u_0(\x, \theta) 
\end{multline*}
or all $(\x, \theta) \in \Omega \times [\underline \theta, \overline \theta].$
Therefore,
\begin{equation}
	\Delta w(\x, \theta) 
	+ \Big[\frac{\Delta u_0(\x, \theta)}{u_0(\x, \theta)} + k^2 c_0(\x)  \Big]w(\x, \theta)
	+ 2 \nabla w(\x, \theta) \cdot \frac{\nabla u_0(\x, \theta)}{u_0(\x, \theta)} = p(\x)
	\label{3.3}
\end{equation}
for all $(\x, \theta) \in \Omega \times [\underline \theta, \overline \theta].$
To eliminate the unknown $p$, we differentiate \eqref{3.3} with respect to $\theta$.
We obtain
\begin{multline}
	\Delta \partial_\theta w(\x, \theta) 
	+ \partial_\theta\Big[\frac{\Delta u_0(\x, \theta)}{u_0(\x, \theta)} + k^2 c_0(\x)  \Big]w(\x, \theta)
	\\
	+\Big[\frac{\Delta u_0(\x, \theta)}{u_0(\x, \theta)} + k^2 c_0(\x)  \Big]\partial_\theta w(\x, \theta)
	+ 2 \nabla \partial_\theta w(\x, \theta) \cdot \frac{\nabla u_0(\x, \theta)}{u_0(\x, \theta)}
	\\
	+ 2 \nabla w(\x, \theta) \cdot  \partial_\theta \Big[\frac{\nabla u_0(\x, \theta)}{u_0(\x, \theta)}\Big]
	 = 0
	\label{3.4}
\end{multline}
for all $(\x, \theta) \in \Omega \times [\underline \theta, \overline \theta].$
Solving the third order equation \eqref{3.4} is extremely challenging. 
An analytic and numerical approach to solve it is not developed yet. 
Suggested by \cite{NguyenLiKlibanov:2019}, we only consider an approximation context due to a truncation of the Fourier series of the solution $w(\x, \theta)$ with respect to the special basis.
This basis was first introduced in \cite{Klibanov:jiip2017}.
We briefly recall  the construction of this basis.
For each $m\geq 1$, define $\phi _{m}(\theta)=(\theta-\theta_{0})^{m-1}\exp (\theta-\theta_{0})$
where $\theta_{0}=(\underline{\theta}+\overline{\theta})/2$. The sequence $\{\phi
_{m}\}_{m=1}^{\infty }$ is complete in $L^{2}(\underline{\theta},\overline{\theta})$.
Applying the Gram-Schmidt orthonormalization procedure to the sequence $%
\{\phi _{m}\}_{m=1}^{\infty }$, we obtain an orthonormal basis in $L^{2}(%
\underline{\theta},\overline{\theta}),$ denoted by $\{\Psi _{m}\}_{m=1}^{\infty }$. It
is not hard to verify that for each $m,$ the function $\Psi _{m}(\theta)$ has the
form
$
\Psi _{m}(\theta)=P_{m-1}(\theta-\theta_{0})\exp (\theta-\theta_{0}),
$
where $P_{m-1}$ is a polynomial of the degree $(m-1)$. The
following result plays an important role in our analysis.

\begin{proposition}[see \protect\cite{Klibanov:jiip2017}]
For $m, n\geq 1$, we have
\begin{equation}
s_{mn}=\int_{\underline{\theta}}^{\overline{\theta}}\Psi _{m}(\theta)\Psi _{n}^{\prime
}(\theta)d\theta=\left\{
\begin{array}{ll}
1 & \mbox{if }n=m, \\
0 & \mbox{if }n<m.%
\end{array}%
\right.  \label{1}
\end{equation}%
Consequently, let $N>1$ be an integer. Then, the $N\times N$ matrix
\begin{equation}
S=(s_{mn})_{m,n=1}^{N}  \label{matrix D}
\end{equation}
has determinant $1$ and is invertible. 
\label{prop MK}
\end{proposition}
%

Recall the function $w$ defined in \eqref{w}.
For each $(\x, \theta) \in \Omega \times [\underline \theta, \overline \theta]$, we write
\begin{equation}
	w(\x, \theta) = \sum_{n = 1}^\infty w_n(\x) \Psi_n(\theta)
	\simeq \sum_{n = 1}^N w_n(\x) \Psi_n(\theta)
	\label{3.7}
\end{equation}
where
\begin{equation}
	w_n(\x) = \int_{\underline \theta}^{\overline \theta} w(\x, \theta) \Psi_n(\theta)d\theta,
	\quad n = 1, 2, \dots
	\label{3.8}
\end{equation}
\begin{remark}[The choice of $N$]
The cut-off number $N$ in \eqref{3.7} is chosen numerically as follows. 
For each $N$, we define 
\[
	\varphi(N) = \Big\|w(\x, \theta) - \sum_{n = 1}^N w_n(\x) \Psi_n(\theta)\Big\|_{L^{\infty}(\Gamma \times [\underline \theta, \overline \theta])}
\]
where $\Gamma$ is a subset of $\partial \Omega.$
Since the function $u$ is given on $\partial \Omega \times [\underline \theta, \overline \theta])$, using \eqref{w} and \eqref{3.8}, we can compute $\varphi(N)$ directly. We then choose $N$ such that $\varphi(N)$ is sufficiently small. In our computation, $N = 35$. With this choice of $N$, we have $\varphi(N) < 5\times10^{-3}.$
See Figure \ref{figchooseN} for an illustration.
\label{rem_N}
\end{remark}

From now on, we assume that the approximation in \eqref{3.7} is valid.
Plugging the truncating formula \eqref{3.7} into \eqref{3.4}, we obtain
\begin{multline}
	\sum_{n = 1}^N \Delta  w_n(\x) \Psi_n'(\theta) 
	+ \partial_\theta\Big[\frac{\Delta u_0(\x, \theta)}{u_0(\x, \theta)} + k^2 c_0(\x)  \Big] \sum_{n = 1}^N w_n(\x) \Psi_n(\theta)
	\\
	+\Big[\frac{\Delta u_0(\x, \theta)}{u_0(\x, \theta)} + k^2 c_0(\x)  \Big]\sum_{n = 1}^N w_n(\x) \Psi_n'(\theta)
	+ 2 \sum_{n = 1}^N \nabla w_n(\x) \Psi_n'(\theta) \cdot \frac{\nabla u_0(\x, \theta)}{u_0(\x, \theta)}
	\\
	+ 2 \sum_{n = 1}^N \nabla w_n(\x) \Psi_n(\theta) \cdot  \partial_\theta \Big[\frac{\nabla u_0(\x, \theta)}{u_0(\x, \theta)}\Big]
	 = 0
	\label{3.9}
\end{multline}
for all $(\x, \theta) \in \Omega \times [\underline \theta, \overline \theta]$.
For each $m \in \{1, \dots, N\},$ we multiply $\Psi_m(\theta)$ to both sides of \eqref{3.9}. Then, we integrate the resulting equation with respect to $\theta$ on $[\underline \theta, \overline \theta]$.
We get
\begin{equation}
	\sum_{n = 1}^N s_{mn}\Delta  w_n(\x)  
	+ \sum_{n = 1}^N a_{mn}(\x)  w_n(\x)
	+ \sum_{n = 1}^N {\bf b}_{mn}(\x) \cdot \nabla  w_n(\x)  = 0
	\label{3.10}
\end{equation}
for all $\x \in \Omega$.
Here,
\begin{align*}
	s_{mn} &= \int_{\underline \theta}^{\overline \theta} \Psi_n'(\theta) \Psi_m(\theta)d\theta,	\\
a_{mn}(\x) &= \int_{\underline \theta}^{\overline \theta} \Big[
		\partial_\theta\Big(\frac{\Delta u_0(\x, \theta)}{u_0(\x, \theta)} + k^2 c_0(\x)  \Big)\Psi_n(\theta)
		\\
		&\hspace{3cm}+ \Big(\frac{\Delta u_0(\x, \theta)}{u_0(\x, \theta)} + k^2 c_0(\x)  \Big)\Psi_n'(\theta)
	\Big]\Psi_m(\theta)d\theta,
	\\
	{\bf b}_{mn}(\x) &= 2 \int_{\underline \theta}^{\overline \theta}  \Big[ \frac{\nabla u_0(\x, \theta)}{u_0(\x, \theta)} \Psi_n'(\theta) +   \partial_\theta \Big( \frac{\nabla u_0(\x, \theta)}{u_0(\x, \theta)} \Big)\Psi_n(\theta)\Big] \Psi_m(\theta)d\theta
\end{align*}
for all $m, n \in \{1, \dots, N\}.$
Equation \eqref{3.10}, with $m \in \{1, \dots, N\}$, forms a linear system of linear elliptic equations for the vector $W = (w_n)_{n = 1}^N.$

We next derive the boundary conditions for $W$.
Due to \eqref{1.2} and \eqref{3.8}, for all $\x \in \partial \Omega$ and $n \in \{1, \dots, N\},$ we have
\begin{equation}
	w_n(\x) = \int_{\underline \theta}^{\overline \theta} \frac{u(\x, \theta)}{-k^2u_0(\x, \theta)} \Psi_n(\theta) d\theta = 
	\int_{\underline \theta}^{\overline \theta} \frac{f(\x, \theta)}{-k^2u_0(\x, \theta)} \Psi_n(\theta) d\theta
	\label{3.11}
\end{equation}
and
\begin{align}
	\partial_{\nu} w_n(\x)  &= \int_{\underline \theta}^{\overline \theta} \partial_{\nu}\Big(\frac{u(\x, \theta)}{-k^2 u_0(\x, \theta)}  \Big) \Psi_n(\theta)d\theta \notag \\
	&= \int_{\underline \theta}^{\overline \theta} \frac{\partial_{\nu} u(\x, \theta) u_0(\x, \theta) - u(\x, \theta)\partial_{\nu}u_0(\x, \theta) }{-k^2 u_0^2(\x, \theta)}\Psi_n(\theta) d\theta \notag\\
	&= \int_{\underline \theta}^{\overline \theta} \frac{g(\x, \theta) u_0(\x, \theta) - f(\x, \theta)\partial_{\nu}u_0(\x, \theta) }{-k^2u_0^2(\x, \theta)} \Psi_n(\theta) d\theta. \label{3.12}
\end{align}
Due to \eqref{3.10}, \eqref{3.11} and \eqref{3.12},
the vector $W$ satisfies
\begin{equation}
\left\{
\begin{array}{ll}
	 \Delta W(\x) + S^{-1}A(\x) W(\x) + S^{-1}{\bf B}(\x) \cdot \nabla W(\x) = 0 & \x \in \Omega,
	\\
	 W(\x)  = F(\x) & \x \in \partial \Omega,\\
	 \partial_{\nu} W(\x) = G(\x)
	 & \x \in \partial \Omega
\end{array}
\right.
\label{3.13}
\end{equation}
where $S = (s_{mn})_{m, n = 1}^N$, $A(\x) = (a_{mn}(\x))_{m, n = 1}^N$, ${\bf B}(\x) = ({\bf b}_{mn}(\x))_{m, n = 1}^N$ and
\begin{align}
	F(\x) &= \displaystyle \Big(\int_{\underline \theta}^{\overline \theta} \frac{f(\x, \theta)}{-k^2 u_0(\x, \theta)} \Psi_n(\theta) d\theta\Big)_{n = 1}^N, \label{3.14}\\
	G(\x) &= \displaystyle \Big(\int_{\underline \theta}^{\overline \theta} \frac{g(\x, \theta) u_0(\x, \theta) - f(\x, \theta)\partial_{\nu}u_0(\x, k) }{-k^2 u_0^2(\x, \theta)} \Psi_n(\theta) d\theta\Big)_{n = 1}^N. \label{3.15}
\end{align}
The invertibility of $S$ is an important property of the basis $\{\Psi_m\}_{m \geq 1}.$ See Proposition \ref{prop MK} and the proof in \cite{Klibanov:jiip2017}.
Using other basis, for e.g the popular trigonometric basis, is not suitable because the corresponding $S$ is not invertible. 
Solving Problem \ref{isp} becomes the problem of finding a function $W$ satisfying \eqref{3.13}.

\begin{remark}
	Using the truncation in \eqref{3.7} to derive \eqref{3.13} is inspired by the Garlekin approximation in the frequency domain. 
	Studying the behavior of this approximation context as $N \to \infty$ is extremely challenging. 
This is out of the scope of the paper.
Although \eqref{3.13} is not exact, the approximation is good enough for us to obtain out of expectation numerical results. 
This phenomenon is true for a long list of inverse problems. Here is an incomplete list of inverse problems that were solved numerically by using similar truncation approaches.
We refer the reader to \cite{KlibanovNguyen:ip2019} for the $X$-ray tomography problem, 
\cite{KlibanovAlexeyNguyen:SISC2019} for an inverse source problem for the full radiative transfer equation, \cite{VoKlibanovNguyen:IP2020, KhoaKlibanovLoc:SIAMImaging2020,  LeNguyen:preprint2021} for the inverse scattering problem,
\cite{Nguyen:CAMWA2020, Nguyens:jiip2020} for a coefficient inverse problem for parabolic equations, \cite{LeNguyen:2020} for an inverse source problem for nonlinear parabolic equations; \cite{LeNguyenNguyenPowell:JOSC2021} for an inverse source problem for hyperbolic equations. Especially, we refer the reader to \cite{NguyenLiKlibanov:2019} for an algorithm to solve an inverse source problem that is a particular case of Problem \ref{isp}. The inverse source problem in \cite{NguyenLiKlibanov:2019} only for the case when $u_0$ does  only on the wave number $k$ while in this paper, $u_0$ is allowed to depend on the spatial variable $\x$. 
\end{remark}

Solving Problem \ref{isp} becomes finding a vector valued function $W$ that satisfies \eqref{3.13}. Since \eqref{3.13} is over-determined, we apply the quasi-reversibility method to solve it.
It is worth mentioning that the quasi-reversibility method was first introduced in \cite{LattesLions:e1969}. Then, it was used very often in solving over-determined boundary value problems, see e.g., \cite{Becacheelal:AIMS2015, Bourgeois:ip2006, Klibanov:jiipp2013, LeNguyenNguyenPowell:JOSC2021, LocNguyen:ip2019, NguyenLiKlibanov:2019}.
Having the solution $W = (w_1, \dots, w_N) \in H^2(\Omega)^N$ of \eqref{3.13} in hand, we can compute $w(\x, \theta)$ for all $(\x, \theta) \in \overline \Omega \times [\underline \theta, \overline \theta]$ via \eqref{3.7}. Then, we can compute the source function $p$ by using \eqref{3.3}. 
This procedure is summarized in Algorithm \ref{alg 2}.

\begin{algorithm}[h!]
	\caption{\label{alg 2} A numerical method to solve Problem \ref{isp}}
	\begin{algorithmic}[1]
	\State Define the basis $\{\Psi\}_{n \geq 1}$ (see \cite{Klibanov:jiip2017}) and choose a cut-off number $N$ as in Remark \ref{rem_N}.
	\State\label{state_quasi} Solve \eqref{3.13} by the quasi-reversibility method (see Section \ref{s5}) for a vector value function $W^{\rm comp} \in H^2(\Omega)^N.$
	\State Compute 
	\[
	    w^{\rm comp}(\x, \theta) = \sum_{n = 1}^N w_n(\x) \Psi_n(\theta)
	    \quad \mbox{for all } (\x, \theta) \in \overline \Omega \times [\underline \theta, \overline \theta].
	\]
	\State The source function is given by
	\begin{multline*}
	    p^{\rm comp}(\x) = \frac{-1}{\overline \theta - \underline \theta} \int_{\underline \theta}^{\overline \theta} 
	    \Big[
	        \Delta w^{\rm comp}(\x, \theta) 
	        + \Big[\frac{\Delta u_0(\x, \theta)}{u_0(\x, \theta)} + k^2 c_0(\x)  \Big]w^{\rm comp}(\x, \theta)
	        \\
	        + 2 \nabla w^{\rm comp}(\x, \theta) \cdot \frac{\nabla u_0(\x, \theta)}{u_0(\x, \theta)}
	    \Big]d\theta
	\end{multline*}
	for all $\x \in \Omega.$
 	\end{algorithmic}
\end{algorithm}

The quasi-reversibility method used in Step \ref{state_quasi} of Algorithm \ref{alg 2} will be presented in Section \ref{s5}.

\section{The quasi-reversibility method} \label{s5}

For the convenience in the analysis, we consider in this section the boundary data $W = F$ and $\partial_{\nu} W = G$ respectively as the ``indirect data" of Problem \ref{isp}.
In fact, having the data of the inverse problem, see \eqref{1.2}, in hand, we can use the explicit formulas in \eqref{3.14} and \eqref{3.15} to compute $F$ and $G$.
Let $\delta > 0$ be a noise level.
Let $F^{\delta}$ and $G^{\delta}$ and $F^*$ and $G^*$ be the noisy and noiseless versions of $F$ and $G$ respectively.
The corresponding noiseless versions are denoted by $F^*$ and $G^*$.
By noise, we mean that there exists an ``error" function $\mathcal E \in H^2(\Omega)^N$ such that
\begin{equation}
    \left\{
    \begin{array}{l}
         \mathcal E|_{\partial \Omega} = F^\delta - F^*,\\
         \partial_{\nu} \mathcal E|_{\partial \Omega} = G^{\delta} - G^*,\\
         \|\mathcal E\|_{H^2(\Omega)} \leq \delta.
    \end{array}
    \right.
    \label{errorfn}
\end{equation}
\begin{remark}
    The existence of the error function in \eqref{errorfn} imply that the noise can be regularized; for e.g. \cite[Section 5]{AmmariWenjiaLoc:jde2013} for the noise model such that it can be smooth out. This condition is significant only for the proof of the convergence theorem (see Theorem \ref{thm_main} and its proof in \cite{NguyenLiKlibanov:2019}). However, in numerical study, we do not have to smooth the noise.
   In fact, we compute the noisy data $F^\delta$ and $G^{\delta}$ using formulas \eqref{3.14} and\eqref{3.15} with $f$ and $g$ replaced by $f^{\delta}$ and $f^{\delta}$ respectively. Here,
    \begin{equation}
        f^{\delta} = f^*(1 + \delta \mbox{rand}),
        \quad
        g^{\delta} = g^*(1 + \delta \mbox{rand})
        \label{noiserand}
    \end{equation}
    where $\mbox{rand}$ is a function taking uniformly distributed random numbers in $[-1, 1]$ and $f^*$ and $g^*$ are the noiseless versions of $f$ and $g$ respectively.
    \label{remnoise}
\end{remark}
Define the set of admissible solution
\begin{equation}  
    H = \big\{
        V \in H^2(\Omega): 
        V|_{\partial \Omega} = F^\delta \mbox{and }
        \partial_{\nu} V|_{\partial \Omega} = G^{\delta}
    \big\}.
    \label{H}
\end{equation}
We assume that the set $H \not = \emptyset.$
Due to the presence of noise in the given data, problem \eqref{3.13} is over-determined. It is natural to solve it using the quasi-reversibility method, which is similar to the least-squares optimization together with a Tikhonov regularization term. That means, we minimize the functional
\[
    J(W) = 
    \int_{\Omega} \big|\Delta W + S^{-1} A(\x) W(\x) + S^{-1} {\bf B}(\x) \cdot \nabla W(\x)\big|^2d\x
    + \epsilon \|W\|_{H^2(\Omega)^N}^2
\]
subject to the boundary conditions $W|_{\partial \Omega} = F^\delta$ and $\partial_{\nu} W|_{\partial \Omega} = G^{\delta}$. 
\begin{remark}
    The presence of the Tikhonov regularization term is significant. In the theoretical part, it makes the functional $J$ coercive, which is important for the existence of a minimizer of $J$.
    In practice, we cannot obtain good numerical result without the presence of this Tikhonov regularization term.
    In our numerical study in Section \ref{s5}, we choose $\epsilon = 10^{-5}$. This value of $\epsilon$ was chosen by a trial and error process. 
    We manually try many values of $\epsilon$ for test 1 until we obtain good numerical result. Then, we use this value for all other tests and for several noise levels $\delta \in [5\%, 50\%]$.
\end{remark}

We have the theorem.
\begin{theorem} 
 Let $F^{\delta}$ and $G^{\delta}$ be the noisy boundary data for $W$. Here, $\delta > 0$ is the noise level in the sense of \eqref{errorfn}. Assume that the set $H$ defined in \eqref{H} is nonempty. Then, for any $\epsilon > 0$, the functional $J$ has a unique minimizer in $H$.
    Denote by $W_{\rm min}$ the obtained minimizer.
     Moreover, let $W^*$ be true solution to \eqref{3.13}. Then, the following estimate holds true
    \begin{equation}
        \|W^{\delta} - W^*\|_{H^1(\Omega)^N}^2 \leq C \Big(\delta^2 + \epsilon \|W^*\|_{H^2(\Omega)^N}^2\Big).
    \end{equation} 
 \label{thm_main}    
\end{theorem}

For brevity, we do not present the proof of Theorem \ref{thm_main}. We refer the reader to \cite[Theorem 3.1 and Theorem 5.1]{NguyenLiKlibanov:2019} for the proof of this theorem. We also refer to 
\cite[Proposition 4.1 and Theorem 4.1]{Nguyen:CAMWA2020} for the proof of a similar theorem when $J$ involves the boundary integrals of the data.

\begin{remark}
	The quasi-reversibility method we employ here is based on least squares optimization with a special Tikhonov regularization term.
	In general, to solve the over-determined boundary value problem \eqref{3.13}, one can use the least squares optimization method with many other choices for the regularization term. It is worth to study the convergence of the method with each of such choice.
\end{remark}

\section{Numerical study}\label{sec_num}

For simplicity in implementation, we numerically solve Problem \ref{isp} in 2D.
Let $\Omega$ be the square  $(-1, 1)^2 \subset \R^2$.
Since solving the Helmholtz equation in \eqref{Heqn} on the whole space, to generate the simulated data, is challenging, we solve its approximation on $\Omega$, say
\begin{equation}
	\left\{
		\begin{array}{ll}
			\Delta u(\x, \theta) + k^2 c_0(\x)u(\x, \theta) = -k^2 p(\x) u_0(\x, \theta) &\x \in \Omega,\\
			\partial_{nu} u(\x, \theta) - \ik u(\x, \theta) = 0 &\x \in \partial \Omega.
		\end{array}
	\right.
	\label{Heqn_appro}
\end{equation}
This change is acceptable in the sense that it does not effect the analysis in Section \ref{sec_Ana} because our arguments depend only on the form of the governing partial differential equation while the boundary values of $u$ and the flux $\partial_{\nu} u$ serve as the given data.
Due to the motivation of Problem \ref{isp} in Section \ref{motivation}, we choose $u_0$ as the incident wave $v_0$, defined in \eqref{incwave} where $k = 3\pi.$ For simplicity, we  choose $c_0 \equiv 1$
The range of the angle $[\underline \theta, \overline \theta] = [0, 2\pi].$
We solve \eqref{Heqn_appro} by the finite difference method.
We arrange a uniform $N_x \times N_x \times N_\theta$ grid of points in $\Omega \times [\underline \theta, \overline \theta]$ as
\begin{multline*}
    \Big\{
        (x_i, y_j, \theta_l): x_i = -1 + (i-1)d_x,  y_j = -1 + (j - 1)d_x,
        \\
        \theta_l = \underline \theta + (l - 1)d_\theta, 
        1 \leq i, j \leq N_x
        1 \leq l \leq N_\theta
    \Big\}
\end{multline*}
where $d_x = 2/(N_x - 1)$ and $d_\theta = (\overline \theta - \underline \theta)/(N_\theta - 1)$. In our computational program $N_x = 80$ and $N_{\theta} = 250.$
Since solving \eqref{Heqn_appro} by the finite difference method to compute the solution $u(\x, \theta)$ is standard in the scientific community, we do not describe the procedure here. 
Having $u|_{\partial \Omega}$ and $\partial_{\nu} u$ in hand, we can compute the indirect data  $F|_{\partial \Omega}$ and $G|_{\partial \Omega}$ \eqref{3.14} and \eqref{3.15} respectively. The noisy data are as in Remark \ref{remnoise} and \eqref{noiserand}.

\begin{figure}[h!]
	\subfloat[]{\includegraphics[width=.3\textwidth]{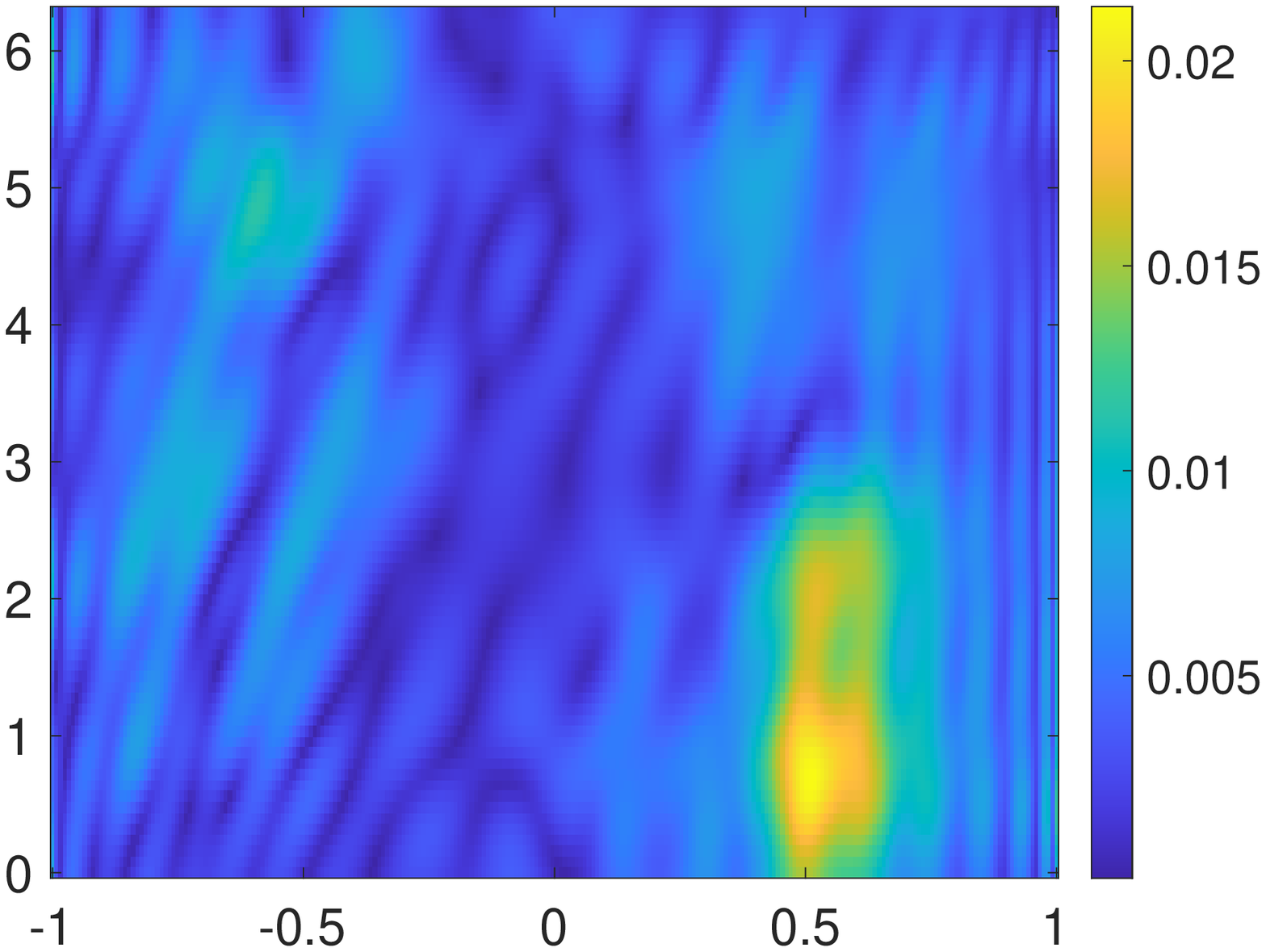}}
	\quad
	\subfloat[]{\includegraphics[width=.3\textwidth]{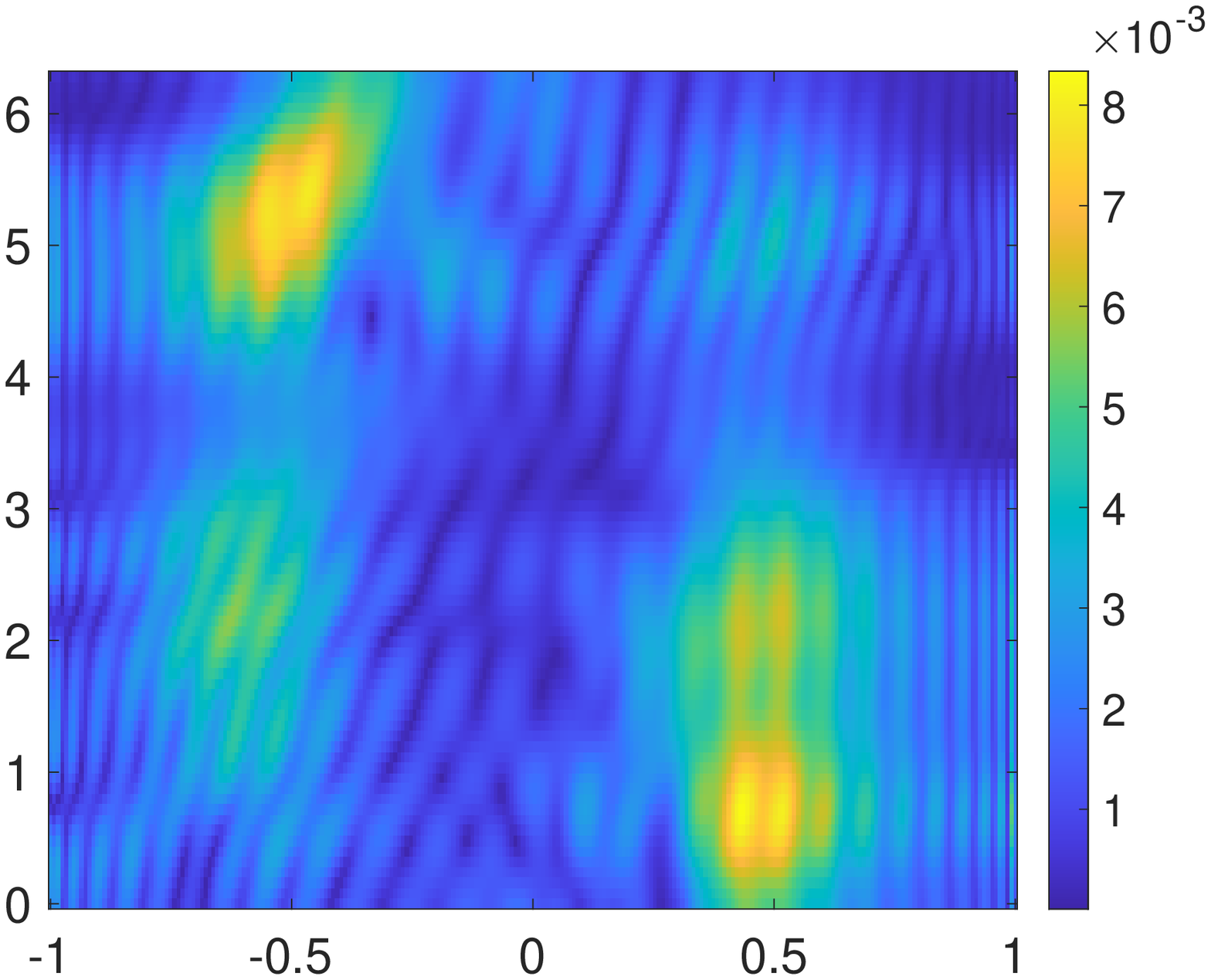}}
	\quad
	\subfloat[]{\includegraphics[width=.3\textwidth]{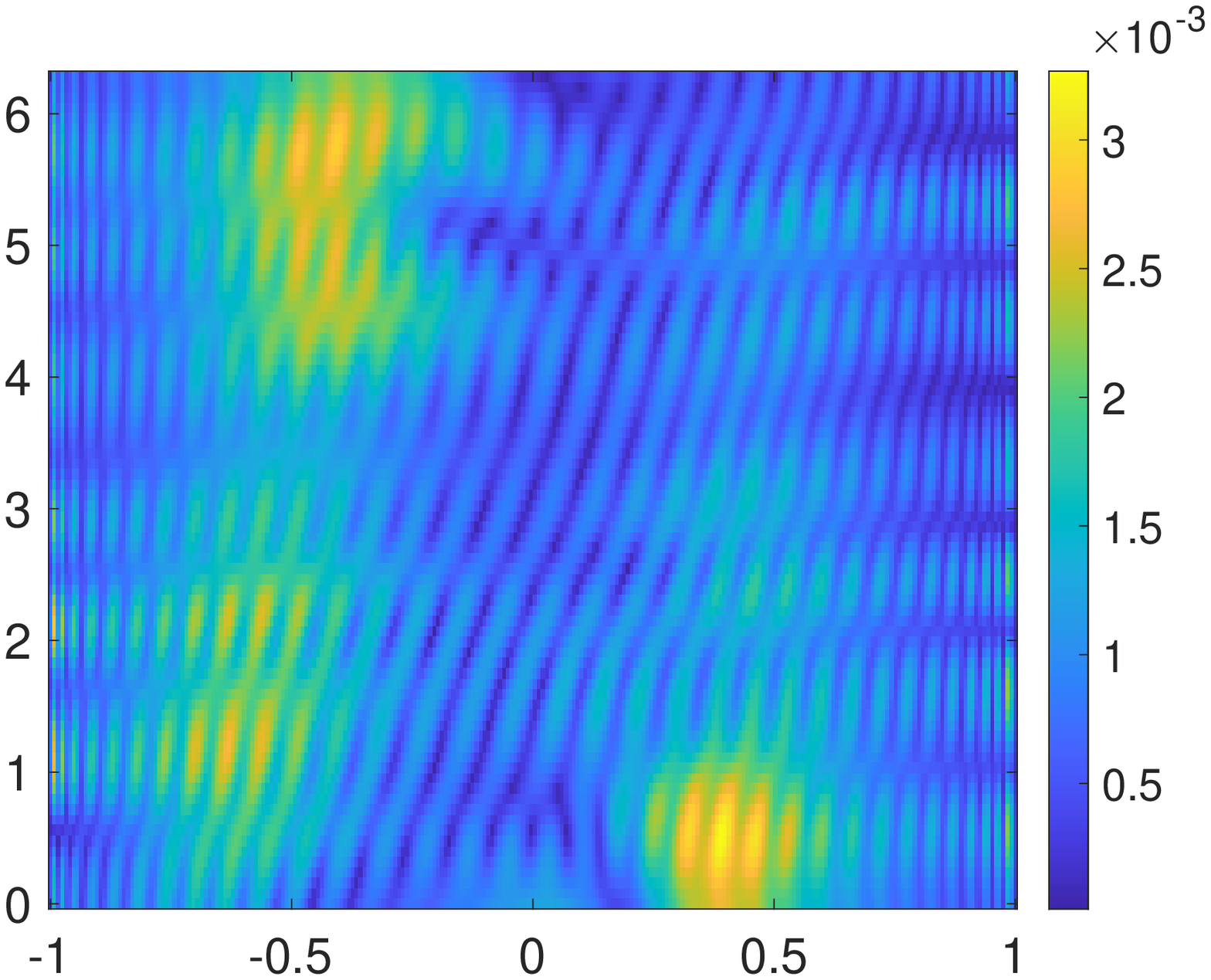}}
	\caption{\label{figchooseN} The graph of the function $\varphi(N) = \|w(\x, \theta) - \sum_{n = 1}^N w_n(\x) \Psi_n(\theta)\|_{L^{\infty}(\Gamma \times [\underline \theta, \overline \theta])}$ where  $\Gamma = \{(x = R, y): |y| < R\}$ is a subset of $\partial \Omega$.
	In (A), $N = 15$ and $\varphi(15) = 21.3 \times 10^{-3}$
	In (B), $N = 25$, $\varphi(N) = 8.34 \times 10^{-3}$
	 In (C), $N = 35$, $\varphi(N) = 3.26 \times 10^{-3}$.
	 The data is taken from Test 4 below. In these graphs, the horizontal axis indicates the range of $y \in [-R, R]$ and the vertical axis represents the angle $\theta \in [0, 2\pi].$}
\end{figure}

We now present an example of the choice of the cut-off number $N$ in Remark \ref{rem_N}. We choose $\Gamma = \big\{(x = 1, y): |y| \leq 1\big\} \subset \partial \Omega$. We numerically examine how the function $\sum_{n = 1}^N w(\x, \theta)$ approximates the function $w(\x, \theta)$ on $\Gamma \times [\underline \theta, \overline \theta]$ by testing their $L^{\infty}$ difference
\[
    \varphi(N) = \Big\| \sum_{n = 1}^N w_n(\x) \Psi_n(\theta) - w(\x)\Big\|_{L^{\infty}(\Gamma \times [\underline \theta, \overline \theta])}.
\]
The graphs of $\varphi(N)$ for $N \in\{ 15, 25, 35\}$ are displayed in Figure \ref{figchooseN}. 
It is evidence that $\varphi(N = 35)$ is sufficiently small. We choose $N = 35$ in all of our numerical tests.

In Step \ref{state_quasi} of Algorithm \ref{alg 2}, rather than solving \eqref{3.13}, we solve the equivalent problem
\begin{equation}
\left\{
\begin{array}{ll}
	 S \Delta W(\x) + A(\x) W(\x) + {\bf B}(\x) \cdot \nabla W(\x) = 0 & \x \in \Omega,
	\\
	 W(\x)  = F(\x) & \x \in \partial \Omega,\\
	 \partial_{\nu} W(\x) = G(\x)
	 & \x \in \partial \Omega.
\end{array}
\right.
\label{4.2}
\end{equation}
We suggest this change to remove the inverse of the matrix $S$ in the system of PDE. The main reason is that the Gram-Schmidt procedure is unstable, so is $S^{-1}.$
The unstability of $S^{-1}$ might lead to some unnecessary difficulties in numerical experiments.
 The implementation of the quasi-reversibility method to solve \eqref{4.2} is very similar to that in \cite{Nguyen:CAMWA2020, Nguyens:jiip2020}. 
 We do not repeat it here.
The implementation of other Steps of Algorithm \ref{alg 2} are straightforward.

We next show four numerical tests.

\begin{figure}[h!]
    \centering
    \subfloat[]{\includegraphics[width=.3\textwidth]{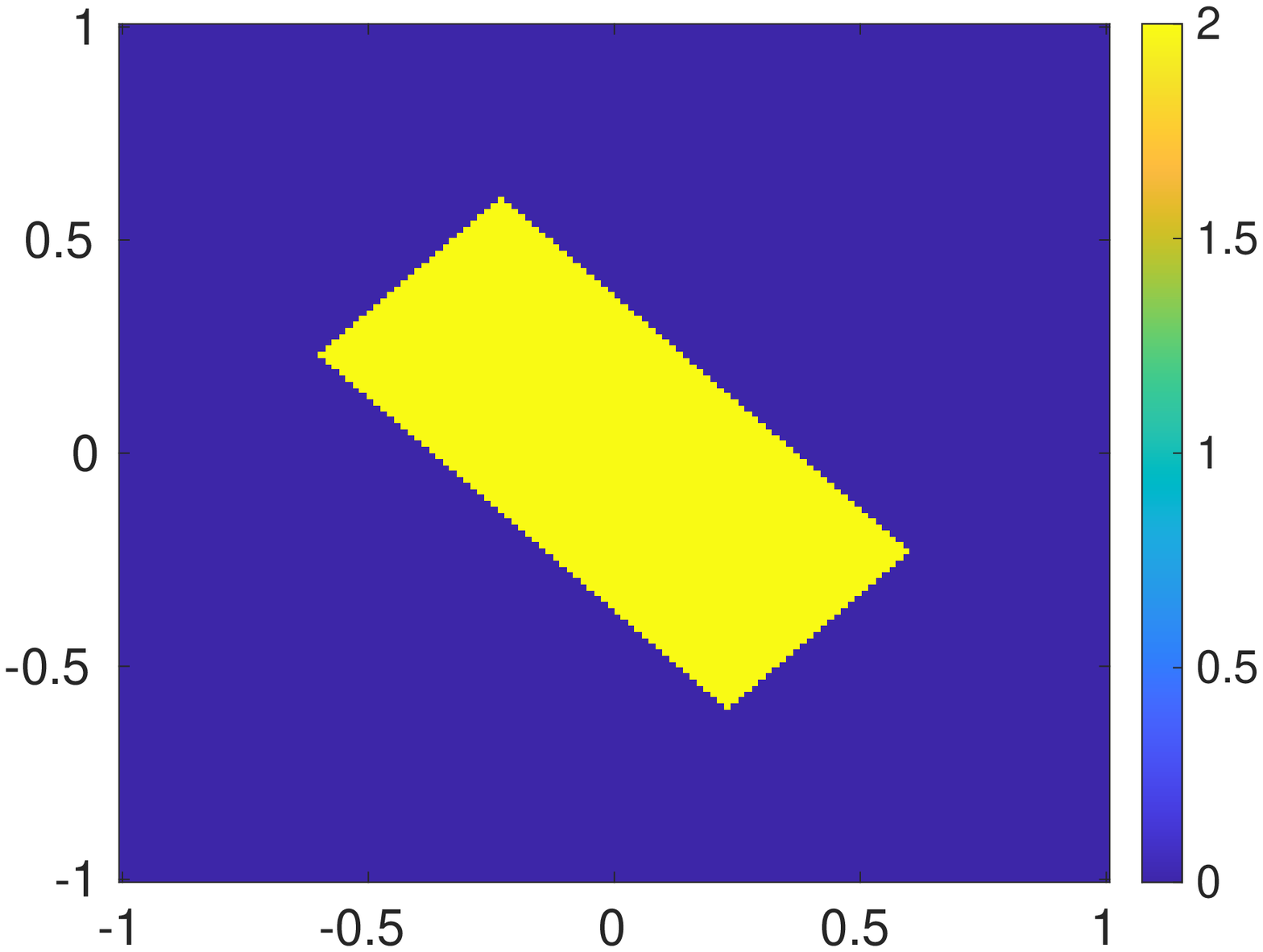}}
    \quad
    \subfloat[]{\includegraphics[width=.3\textwidth]{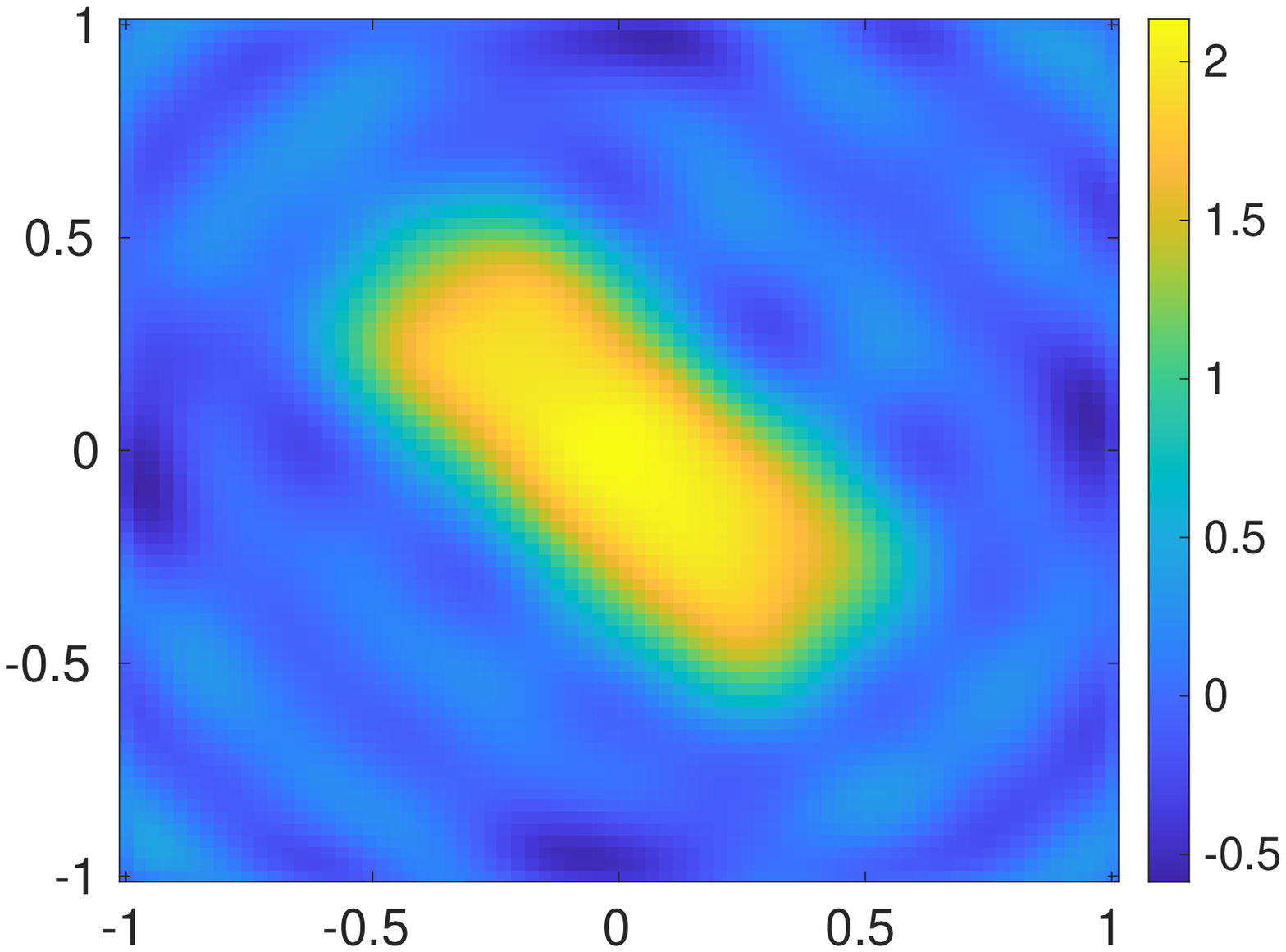}}
    \quad
    \subfloat[]{\includegraphics[width=.3\textwidth]{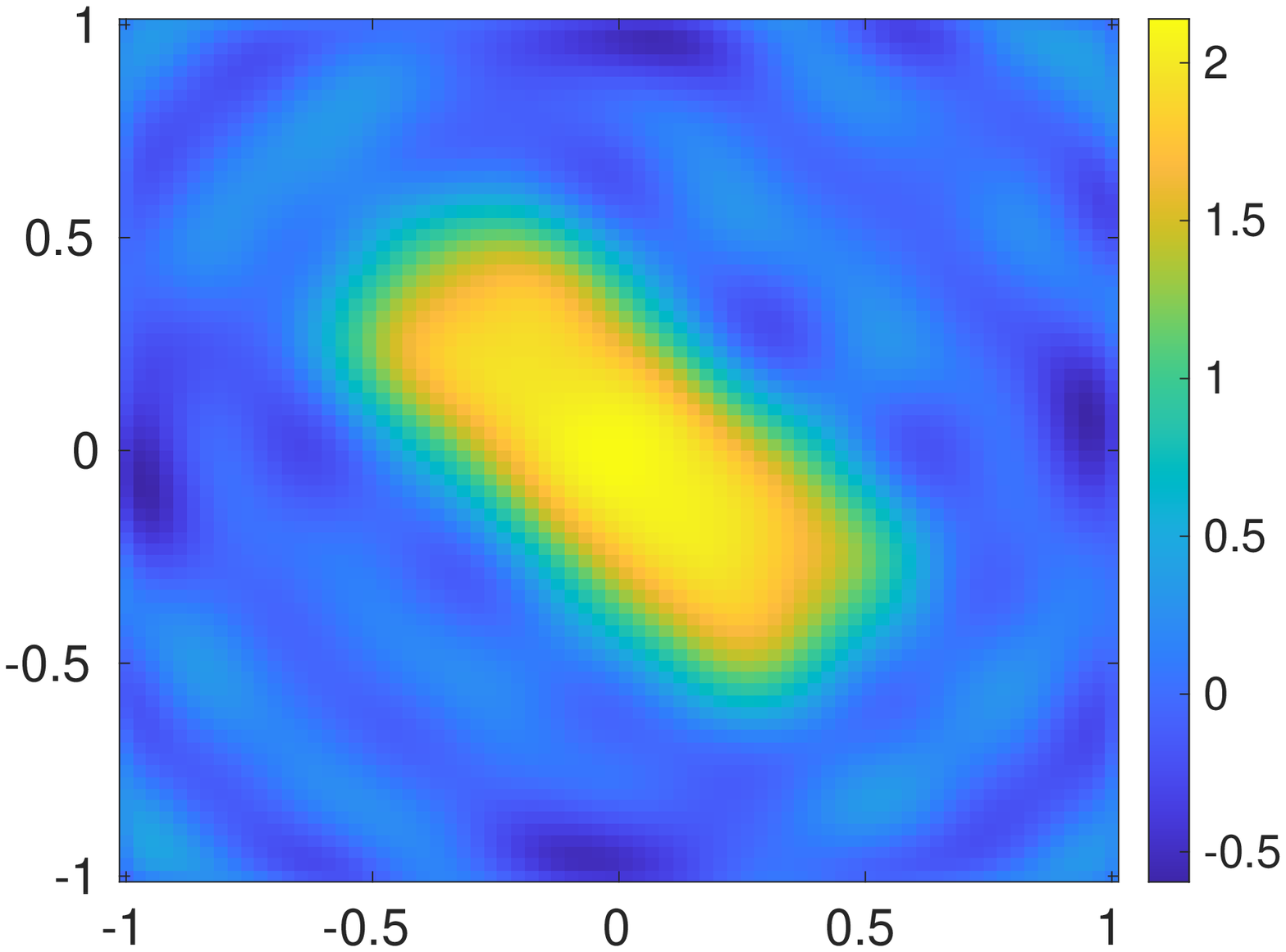}}
    \caption{Test 1. The true and reconstruction of the source functions $p$. (A) The function $p_{\rm true}$. (B) The function $p_{\rm comp}$ computed from data containing $\delta = 5\%$ noise.
    (C) The function $p_{\rm comp}$ computed from data containing $\delta = 10\%$ noise.}
    \label{fig_test1}
\end{figure}

\noindent{\bf Test 1.} We test the case when the support of the true source function takes the form
\[
    p_{\rm true}(x, y) = \left\{
        \begin{array}{ll}
             2&\mbox{if }\max(|x - y|/0.8,|x + y|/0.35)  < 1, \\
             0& \mbox{otherwise}
        \end{array}
    \right.
    \quad \mbox{for all } (x, y) \in \Omega.
\]
The graphs of the true source function and its reconstructions are displayed in Figure \ref{fig_test1}.
 In this test, the support of the true source function $p$ looks like a rectangle. It is evident that this rectangle is detected successfully. 
 On the other hand, the maximum values of the source 
are computed quite accurate. When $\delta = 5\%$, 
$\frac{|\max_{\x \in \overline \Omega} p_{\rm comp} - \max_{\x \in \overline \Omega} p_{\rm true}|}{|\max_{\x \in \overline \Omega} p_{\rm true}|}  = 6.72\%$. 
When $\delta = 10\%$, this error is $6.95\%$.

\begin{figure}[h!]
    \centering
    \subfloat[]{\includegraphics[width=.3\textwidth]{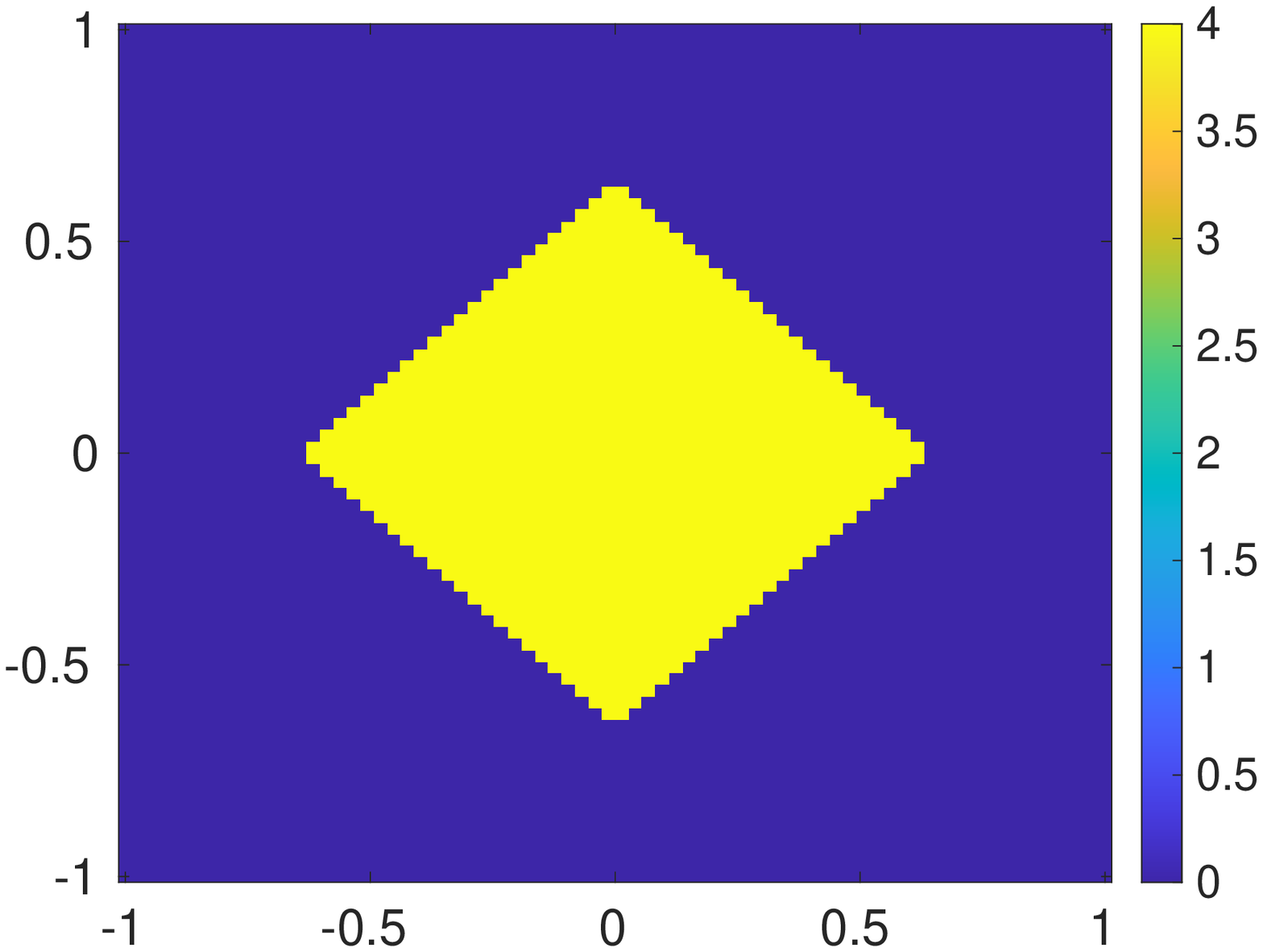}}
    \quad
    \subfloat[]{\includegraphics[width=.3\textwidth]{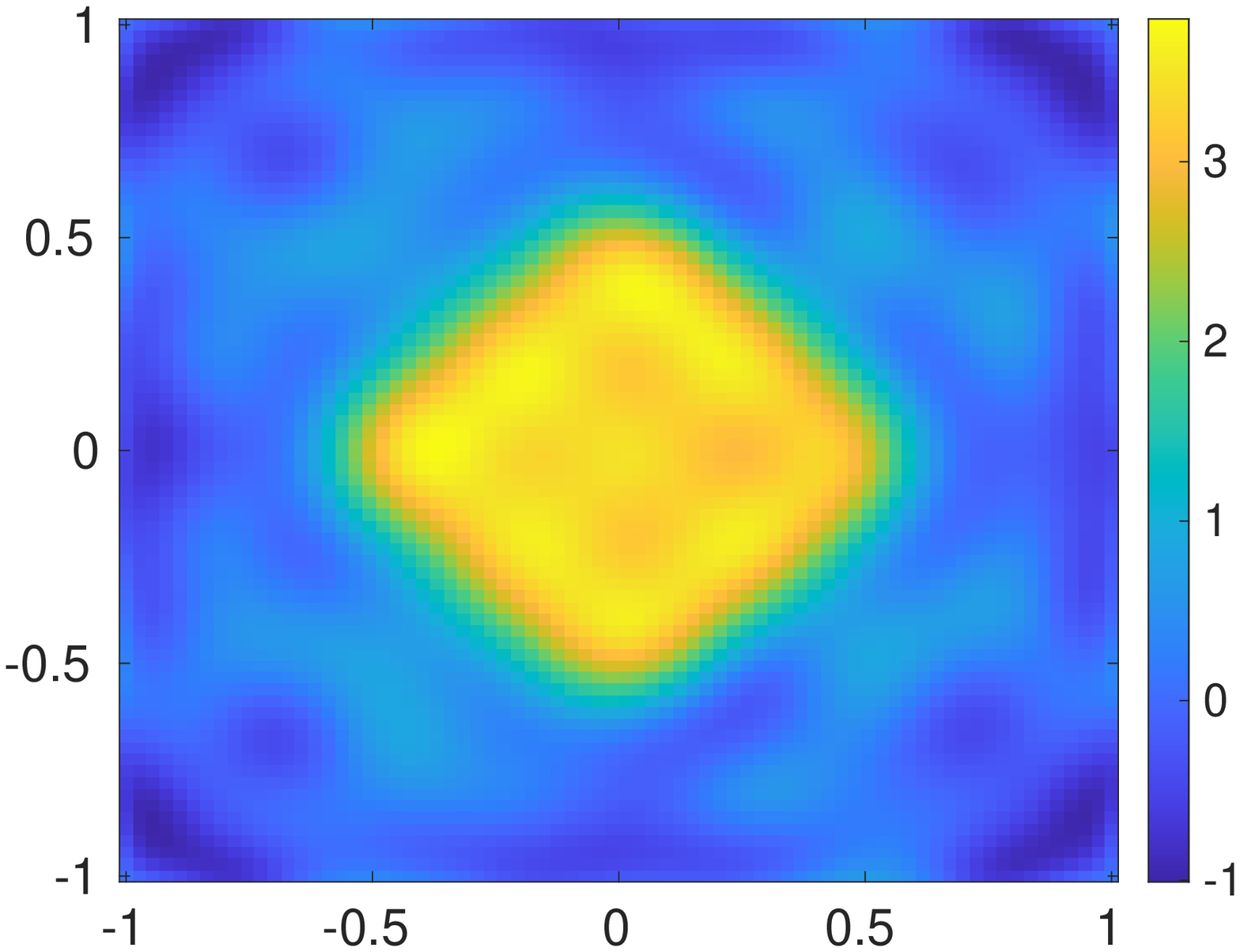}}
    \quad
    \subfloat[]{\includegraphics[width=.3\textwidth]{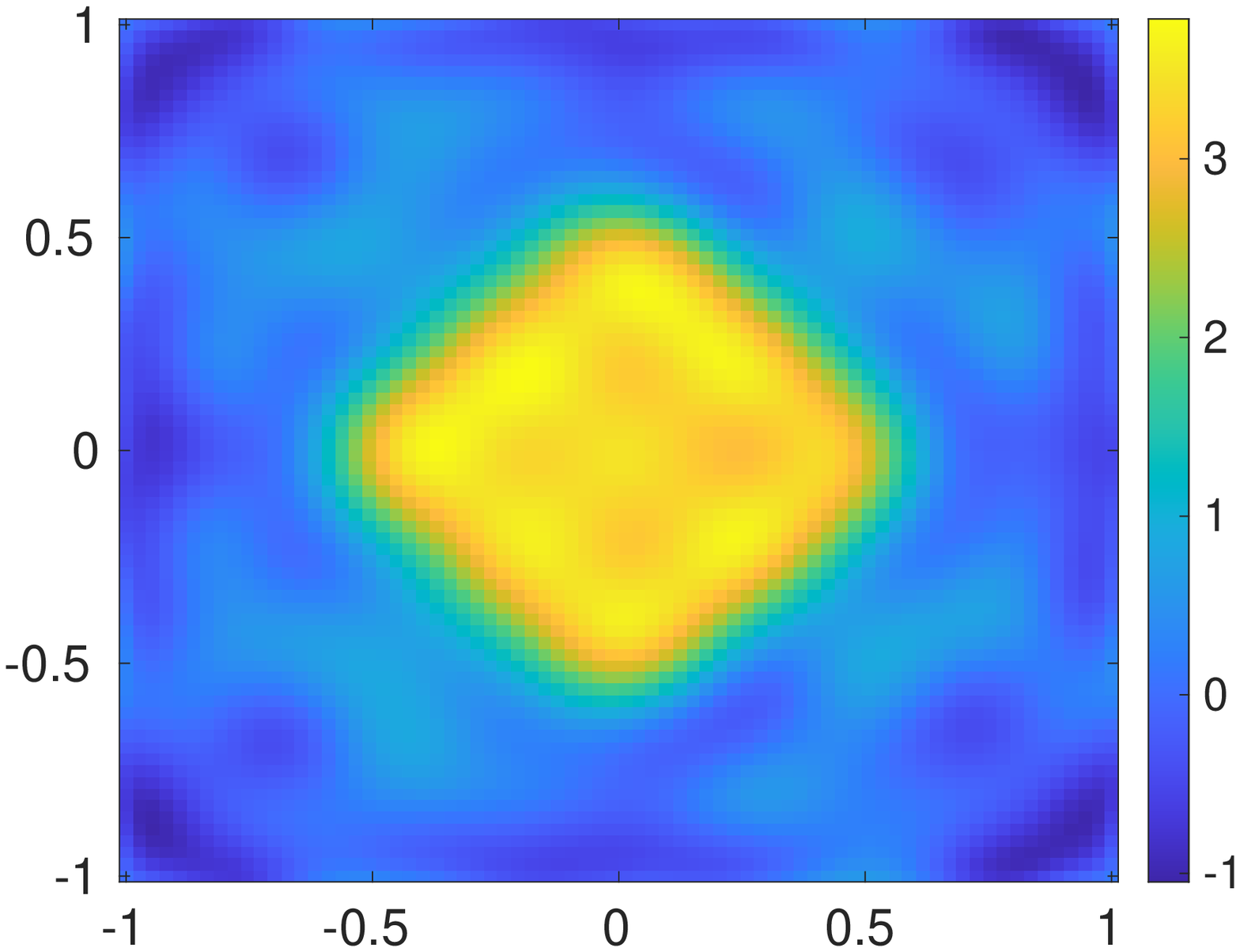}}
    \caption{Test 2. The true and reconstruction of the source functions $p$. (A) The function $p_{\rm true}$. (B) The function $p_{\rm comp}$ computed from data containing $\delta = 10\%$ noise.
    (C) The function $p_{\rm comp}$ computed from data containing $\delta = 30\%$ noise.}
    \label{fig_test2}
\end{figure}

\noindent {\bf Test 2.} We next test the case when the value and the support of the true source are larger than those in Test 1.
The true source function is given by
\[
    p_{\rm true}(x, y) = \left\{
        \begin{array}{ll}
             4 & \mbox{if } \max\{|x + y|, |x - y|\} < 0.6\\
             0 & \mbox{otherwise}
        \end{array}
    \right.
    \quad \mbox{for all } (x, y) \in \Omega.
\]
The graphs of the true source function and its reconstructions are displayed in Figure \ref{fig_test2}.
The support of true source function $p_{\rm true}$ is the square centered at the origin, rotated $45^\circ$. We can see that reconstruction this square is out of expectation. The maximum values the source function inside the square is quite correctly computed. 
When $\delta = 10\%$, $
\frac{|\max_{\x \in \overline \Omega} p_{\rm comp} - \max_{\x \in \overline \Omega} p_{\rm true}|}{|\max_{\x \in \overline \Omega} p_{\rm true}|}= 5.13\%$. 
When $\delta = 30\%$, this error is $5.44\%$.

\begin{figure}[h!]
    \centering
    \subfloat[]{\includegraphics[width=.3\textwidth]{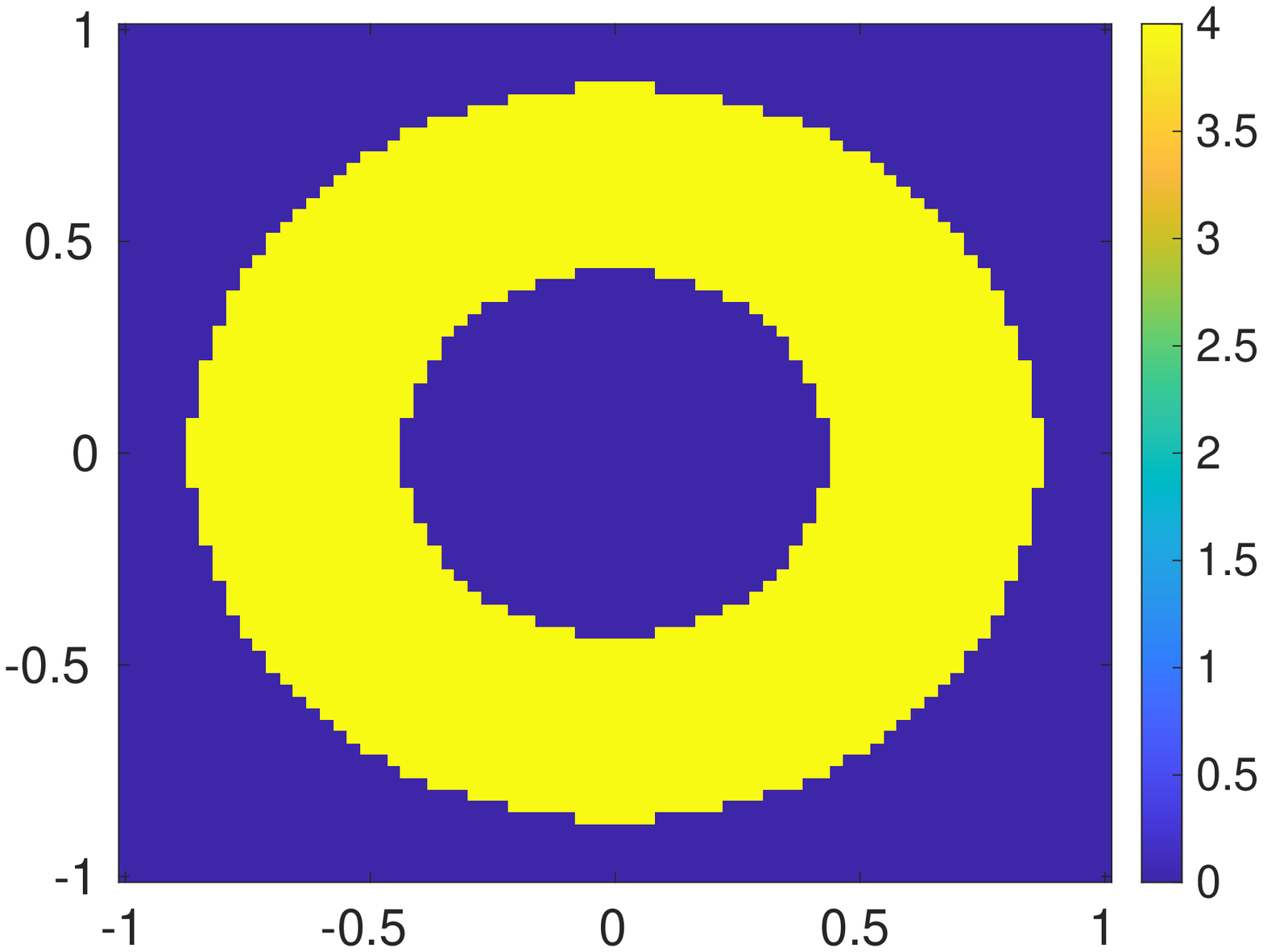}}
    \quad
    \subfloat[]{\includegraphics[width=.3\textwidth]{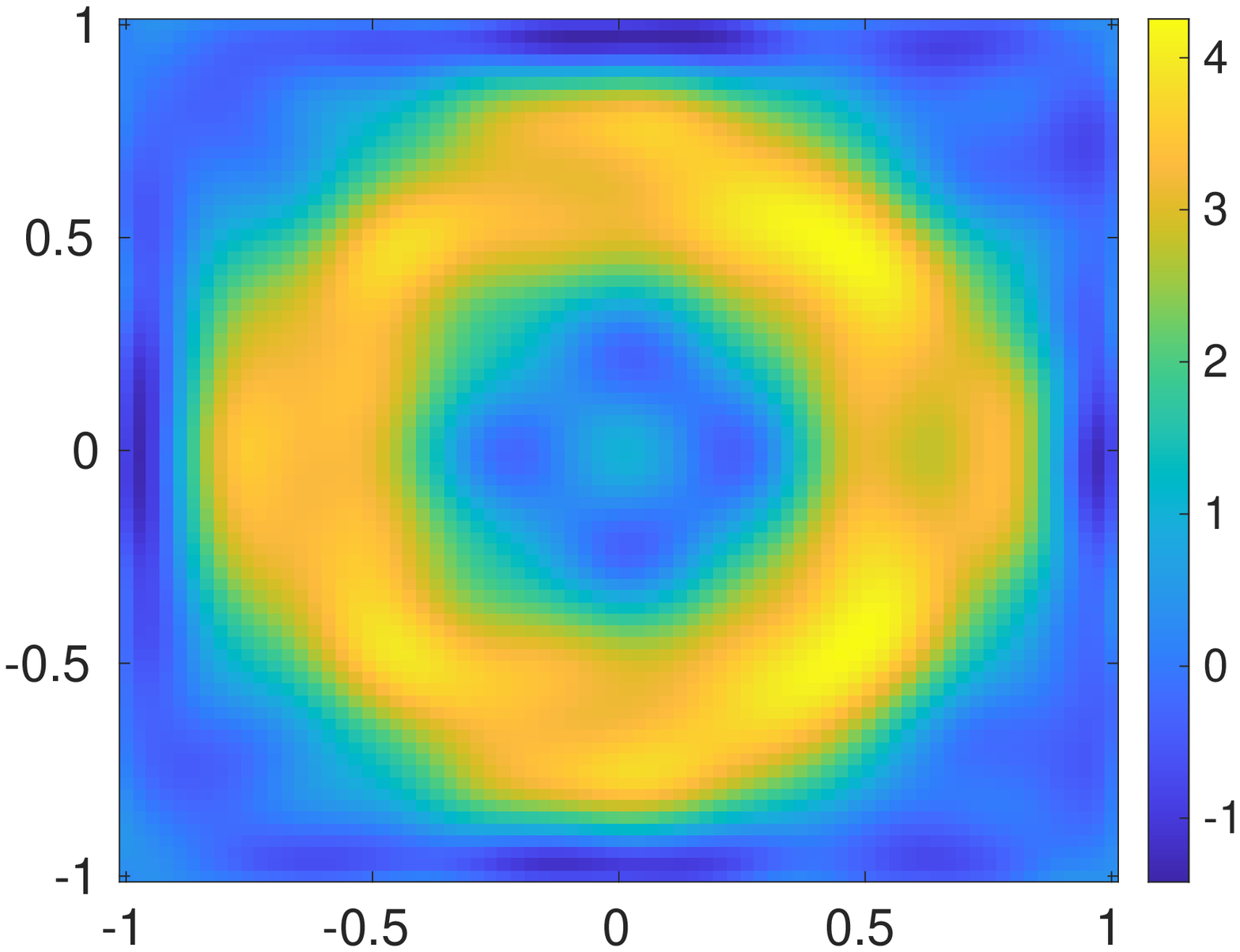}}
    \quad
    \subfloat[]{\includegraphics[width=.3\textwidth]{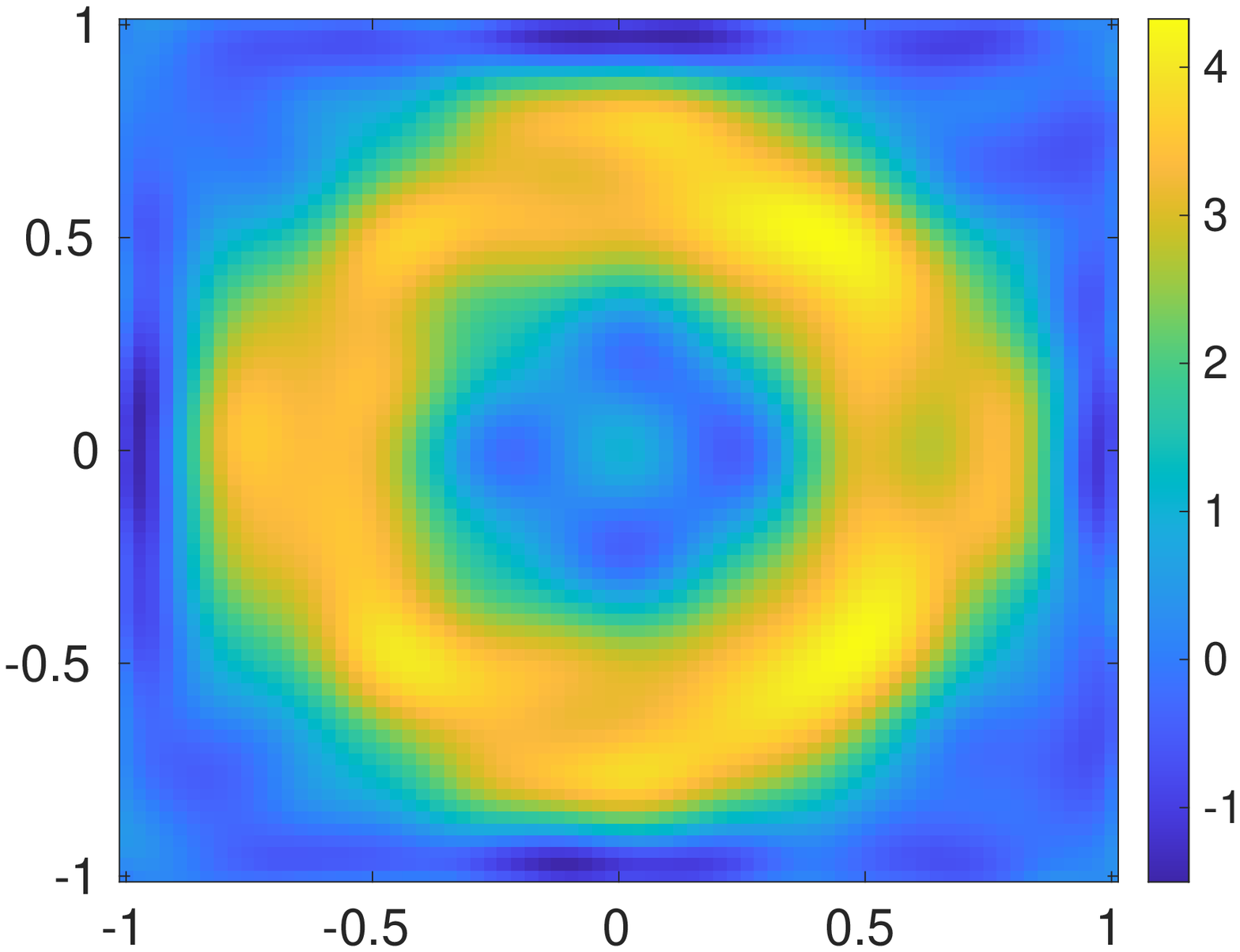}}
    \caption{Test 3. The true and reconstruction of the source functions $p$. (A) The function $p_{\rm true}$. (B) The function $p_{\rm comp}$ computed from data containing $\delta = 5\%$ noise.
    (C) The function $p_{\rm comp}$ computed from data containing $\delta = 50\%$ noise.}
    \label{fig_test3}
\end{figure}

\noindent{\bf Test 3.} We next test the case when the graph of the true source has a void.
The true source function in this test is
\[
    p_{\rm true}(x, y) = \left\{
        \begin{array}{ll}
             4 &  \mbox{if } 0.4^2 < x^2 + y^2 < 0.8^2\\
             0 & \mbox{otherwise}
        \end{array}
    \right.
    \quad \mbox{for all } (x, y) \in \Omega.
\]
The graphs of the true source function and its reconstructions are displayed in Figure \ref{fig_test3}.
The support of true source function $p_{\rm true}$ is the ring centered at the origin with outer radius $0.8$ and inner radius $0.4$. 
This test is interesting due to the presence of the void. 
We can see that reconstruction this ring is acceptable. The maximum values the source function inside the square is quite correctly computed. 
When $\delta = 5\%$, $\frac{|\max_{\x \in \overline \Omega} p_{\rm comp} - \max_{\x \in \overline \Omega} p_{\rm true}|}{|\max_{\x \in \overline \Omega} p_{\rm true}|}  = 6.35\%$. 
When $\delta = 50\%$, this error is $8.13\%$.

\begin{figure}[h!]
    \centering
    \subfloat[]{\includegraphics[width=.3\textwidth]{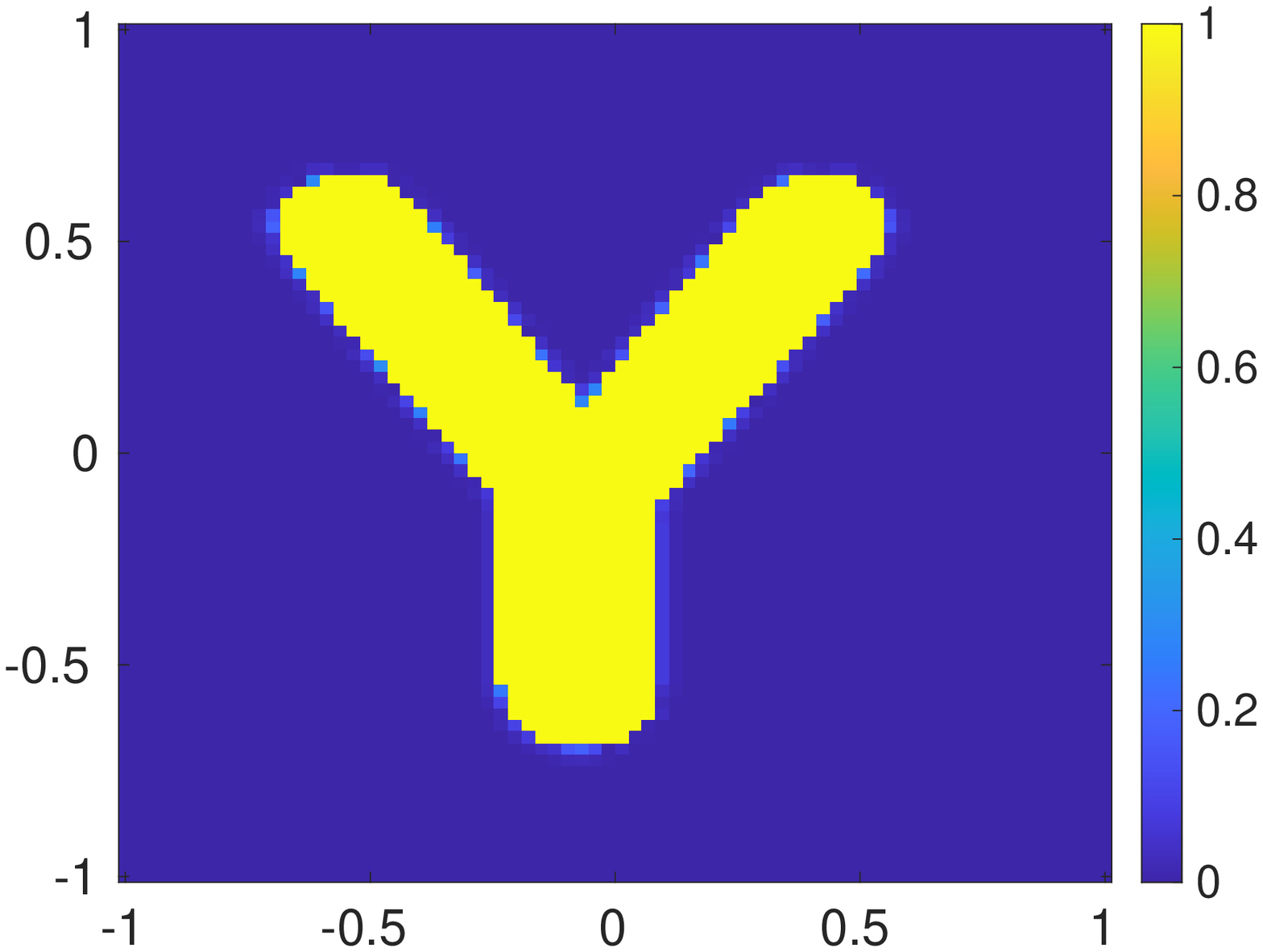}}
    \quad
    \subfloat[]{\includegraphics[width=.3\textwidth]{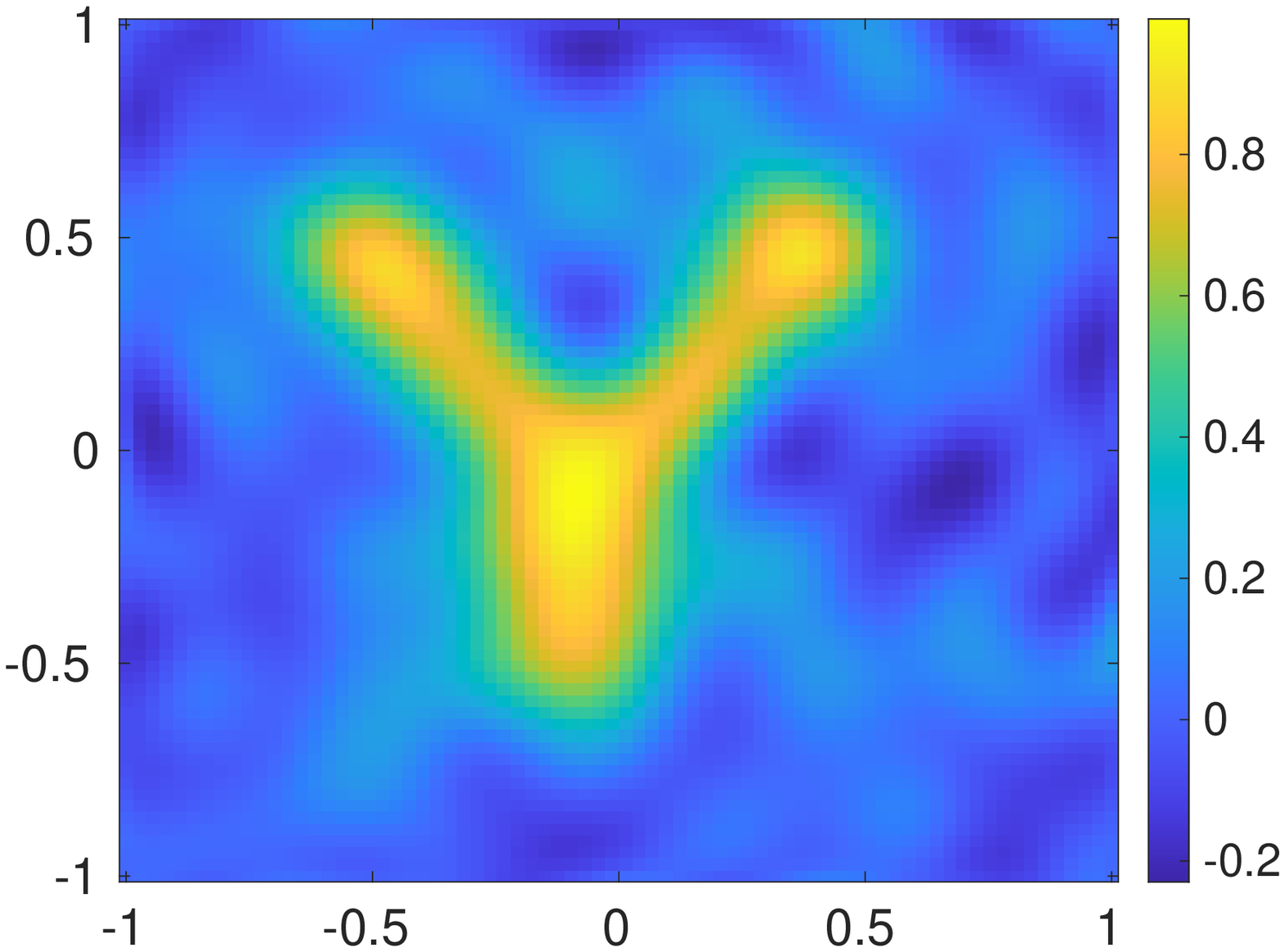}}
    \quad
    \subfloat[]{\includegraphics[width=.3\textwidth]{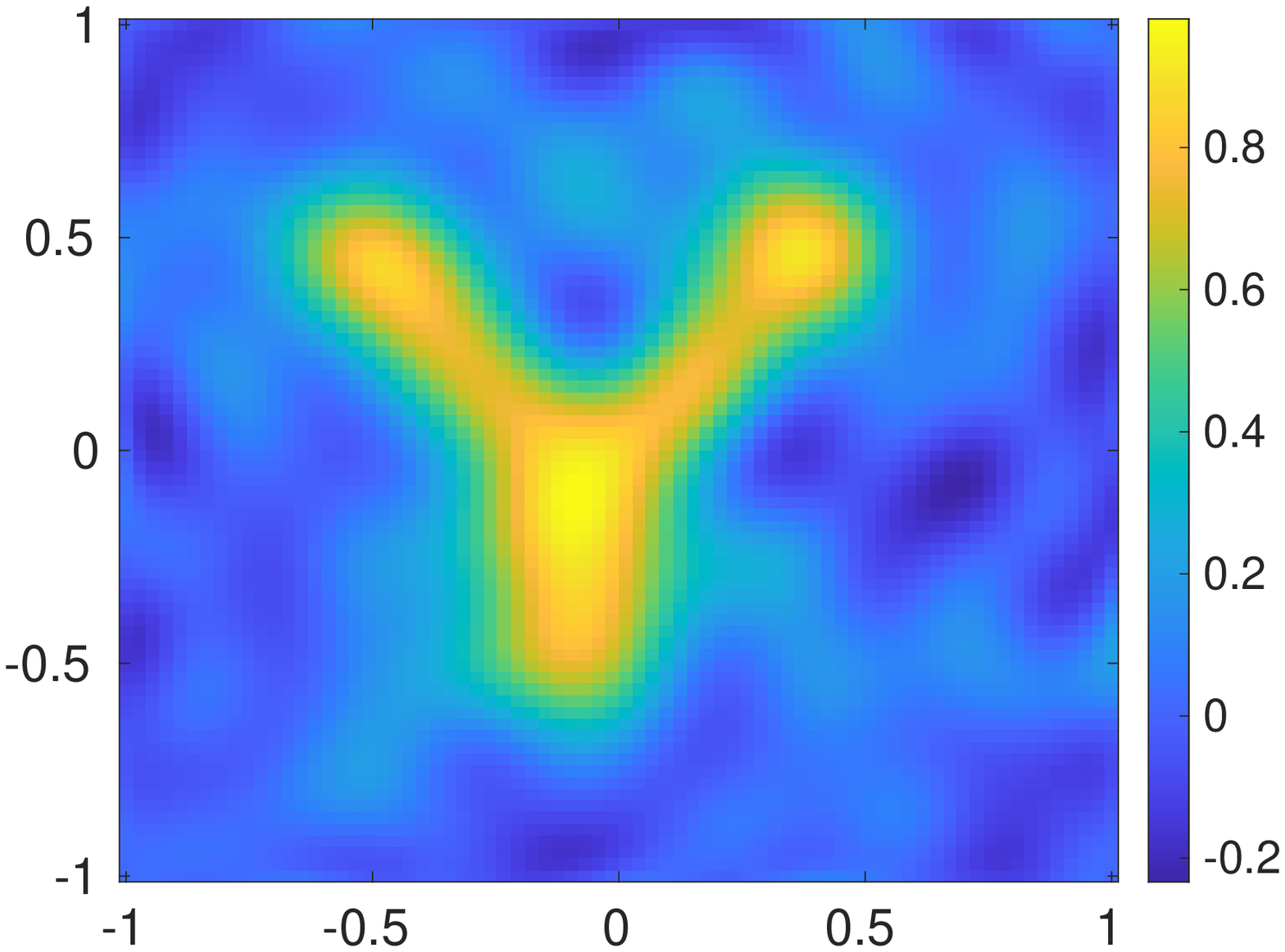}}
    \caption{Test 4. The true and reconstruction of the source functions $p$. (A) The function $p_{\rm true}$. (B) The function $p_{\rm comp}$ computed from data containing $\delta = 10\%$ noise.
    (C) The function $p_{\rm comp}$ computed from data containing $\delta = 50\%$ noise.}
    \label{fig_test4}
\end{figure}

\noindent{\bf Test 4.}
We next consider a more interesting case. The function $p_{\rm true}$ is the characteristic function of the letter $Y$. Although this true source function has complicated structure, we are able to well compute it.
The graphs of the true source function and its reconstructions are displayed in Figure \ref{fig_test4}.
We can see that reconstruction the letter $Y$ is acceptable. The maximum values the source function inside the square is quite correctly computed. 
When $\delta = 10\%$, $\frac{|\max_{\x \in \overline \Omega} p_{\rm comp} - \max_{\x \in \overline \Omega} p_{\rm true}|}{|\max_{\x \in \overline \Omega} p_{\rm true}|}  = 0.83\%$. 
When $\delta = 50\%$, this error is $1.92\%$.

\begin{remark}
These numerical examples numerically show that Algorithm \ref{alg 2} is robust. Especially, it is  stable with respect to the noise. We can obtain good numerical results even when the noise level is up to $50\%$.
\end{remark}

\section{Concluding remarks}\label{sec_concl}
In this paper, we solve an inverse source problem. 
This inverse source problem is the linearization of the highly nonlinear and severely ill-posed inverse scattering problem. 
In order to solve the inverse source problem, in the first step, we derive a system of linear elliptic PDEs in which the source function is absent.
We solve this system by the quasi-reversibility method. 
We choose the quasi-reversibility method because its convergence as noise tends to $0$ was proved.
The efficiency of our method is confirmed by some interesting numerical examples.



\end{document}